\documentclass[12pt,reqno]{article}
\usepackage{amssymb,amscd,amsmath,amsthm,color}
\newtheorem{theorem}{Theorem}[section]
\newtheorem{lemma}[theorem]{Lemma}
\newtheorem{cor}[theorem]{Corollary}

\usepackage{enumerate,amssymb}
\theoremstyle{definition}

\theoremstyle{remark}
\newtheorem{remark}[theorem]{Remark}

\numberwithin{equation}{section}

\def\bC{\mathbb{C}}

\def\bM{\mathbb{M}}
\def\bN{\mathbb{N}}
\def\bR{\mathbb{R}}

\begin{document}
\baselineskip=15pt

\title{Decomposition and partial trace of  positive   matrices with Hermitian blocks }

\author{Jean-Christophe Bourin{\footnote{Supported by ANR 2011-BS01-008-01.}}\, and  Eun-Young Lee{\footnote{Research  supported by  Basic Science Research Program
through the National Research Foundation of Korea (NRF) funded by the Ministry of Education,
Science and Technology (2010-0003520)}}}

\date{ }

\maketitle

\begin{abstract}
\noindent  Let $H=[A_{s,t}]$ be a positive definite matrix written in $\beta\times\beta$ Hermitian blocks and let $\Delta=A_{1,1}+\cdots+A_{\beta,\beta}$ be its partial trace.   Assume that $\beta=2^p$ for some  $p\in\bN$. Then, up to a direct sum operation,   $ H$ is the average of $\beta$ matrices isometrically congruent to $ \Delta$. 
  A few corollaries are given, related to important  inequalities in quantum information theory such as the Nielsen-Kempe separability criterion. 
\end{abstract}

{\small\noindent

\noindent
AMS subjects classification 2010:  15A60, 47A30,  15A42.}

\section{Introduction and a key lemma}

 Positive semi-definite matrices partitioned in two by two blocks   occur as an efficient tool in matrix analysis, sometimes a magic tool ! $-$ according Bhatia's famous book \cite{Bhatia}. These partitions allow to derive a lot of important inequalities and those with Hermitian blocks shed much light on the geometric and harmonic matrix means. Partitions into a larger number of blocks  are naturally involved with tensor products, in the theory of positive linear maps and in their application in quantum physics.

This article deals with positive matrices partitioned in Hermitian blocks.  By using  unitary or isometry congruences, we will  improve some nice majorisations, or norm estimates,  first obtained in the field of quantum information theory.

For  partitioned positive matrices,
 the diagonal blocks  play a quite special role. This is apparent in a rather striking  decomposition  due  to  the  authors \cite{BL1}.

\begin{lemma} \label{BL-lemma} For every matrix in  $\bM_{n+m}^+$ written in blocks, we have a decomposition
\begin{equation*}
\begin{bmatrix} A &X \\
X^* &B\end{bmatrix} = U
\begin{bmatrix} A &0 \\
0 &0\end{bmatrix} U^* +
V\begin{bmatrix} 0 &0 \\
0 &B\end{bmatrix} V^*
\end{equation*}
for some unitaries $U,\,V\in  \bM_{n+m}$.
\end{lemma}

This lemma leads  to  study  partitions via unitary congruences. It  is the key of the subsequent results. A proof and several consequences can be found in \cite{BL1} and \cite{BH1}. 
Of course,  $\bM_n$ is the algebra of $n\times n$  matrices with real or complex entries, and $\bM_n^+$ is the positive  part. That is, $\bM_n$ may stand either for $\bM_n(\bR)$, the matrices with real entries, or for  $\bM_n(\bC)$, those with complex entries. The situation is different
in the next statement, where complex entries seem unavoidable. 

\vskip 10pt
\begin{theorem}\label{thm-four} Given any matrix in $\bM_{2n}^+(\bC)$ written in blocks in $\bM_n(\bC)$ with Hermitian off-diagonal blocks, we have
\begin{equation*}
\begin{bmatrix} A &X \\
X &B\end{bmatrix}= \frac{1}{2}\left\{ U(A+B)U^* +V(A+B)V^*\right\}
\end{equation*} for some isometries $U,V\in\bM_{2n,n}(\bC)$.
\end{theorem}

\vskip 10pt
Here $\bM_{p,q}(\bC)$ denote the space of $p$ rows and $q$ columns matrices with complex entries, and $V\in\bM_{p,q}(\bC)$ is  an isometry if $p\ge q$ and $V^*V=I_q$. Even for a matrix in $\bM_{2n}^+(\bR)$, it seems essential to use isometries with complex entries ! The result, due to  Lin and the authors, is based on Lemma  \ref{BL-lemma}, a proof is in \cite{BLL2} and implicitly in \cite{BLL1}.

There is no evidence whether a positive block-matrix $H$ in $\bM^+_{3n}$, 
$$
H=\begin{bmatrix}
A&X&Y \\ X&B&Z \\Y&Z&C
\end{bmatrix}
$$
with Hermitian off-diagonal blocks $X,Y,Z$, could be decomposed as
$$
H=\frac{1}{3}\left\{ U\Delta U^*+ V\Delta V^* +W\Delta W^*\right\}
$$
where $\Delta=A+B+C$ and $U,V,W$ are isometries. In fact, this would be  surprising. However,  a quite nice decomposition is possible by considering direct sum copies: this provides a substitute to Theorem \ref{thm-four} for partitions into an arbitrary number of blocks; it is the main result of this article. 

These decompositions entail  some nice  inequalities. Lemma \ref{BL-lemma} yields a simple estimate for all symmetric (or unitarily invariant) norms,
\begin{equation}\label{fact-0}
\left\|\begin{bmatrix} A &X \\
X^* &B\end{bmatrix} \right\|\le \| A\|+\| B \|. 
\end{equation}
Recall that a symmetric norm on $\bM_m$ satisfies $\|A\|=\|UA\|=\|AU\|$ for all $A\in\bM_m$ and all unitaries $U\in\bM_m$. This obviously induces a symmetric norm on $\bM_n$, $1\le n\le m$. The most familiar symmetric norms are the Schatten $p$-norms, $1\le p< \infty$,
\begin{equation}\label{Schatt}
\| A\|_p = \{{\mathrm{Tr\,}} (A^*A)^{p/2}\}^{1/p},
\end{equation}
and, with $p\to\infty$,  the operator norm. In general,  the sum of the norms $\|A\| +\| B\|$ can not be replaced in \eqref{fact-0} by the norm of the sum $\| A+B \|$. However,  Theorem \ref{thm-four} implies the following remakable corollary.

\begin{cor} Given any matrix in $\bM_{2n}^+$ written in blocks in $\bM_n$ with Hermitian off-diagonal blocks, we have
\begin{equation*}
\left\|\begin{bmatrix} A &X \\
X &B\end{bmatrix} \right\|\le \| A+B \|
\end{equation*} for all symmetric norms.
\end{cor}

This is the simplest case of Hiroshima's theorem, discussed in the next section.
 There are some positive  matrices in $\bM_6$ partitioned in blocks in $\bM_3$, with  normal off-diagonal blocks $X$, $X^*$, such that
\begin{equation*}
\left\| \begin{pmatrix} A& X\\ X^*&B \end{pmatrix} \right\|_{\infty} > \| A+B\|_{\infty}.
\end{equation*}
Hence the assumptions  are rather optimal. 

In Section 2, we state our decomposition  and derive several inequalities, most of them related to Hiroshima's theorem. The proof of the decomposition is given in Section 3. A discussion of previous results for small partitions and some remarks related to quantum information are given in the last section.

\section{Direct sum and partial trace}

A typical example of positive matrices written in blocks are formed by tensor products. Indeed, 
 the tensor product $A\otimes B$ of $A\in\bM_{\beta}$ with $B\in\bM_n$ can be identified with an element of $\bM_{\beta}(\bM_n)=\bM_{\beta n}$. Starting with positive matrices in  $\bM_{\beta}^+$ and $\bM_n^+$ we then get a matrix in $\bM_{\beta n}^+$ partitioned in blocks in $\bM_n$. In quantum physics, sums of tensor products of positive semi-definite (with trace one) occur as so-called separable states. In this setting of tensor products, the sum of the diagonal block is called the partial trace (with respect to $\bM_{\beta}$). We will use this terminology.

\vskip 10pt
\begin{theorem}\label{thm-direct} Let $H=[A_{s,t}]\in \bM_{\beta n}^+$ be written in $\beta\times \beta$  Hermitian blocks in $\bM_n$ and let $\Delta=\sum_{s=1}^{\beta}{A_{s,s}}$ be its partial trace.    If $\beta$ is dyadic, then, with  $m=2^{\beta}$, we have
\begin{equation*}
\oplus^{m} H =\frac{1}{\beta} \sum_{k=1}^{\beta} V_k\left(\oplus^{m}\Delta\right) V_k^*
\end{equation*} 
where   $\{V_k\}_{k=1}^{\beta}$  is a family of isometries in $\bM_{m\beta n,mn}$.
\end{theorem}

\vskip 10pt
Here the spaces $\bM_n$ and $\bM_{p,q}$ denote either  the real or  complex spaces of matrices. By a dyadic number $\beta$, we mean  $\beta=2^p$ for some $p\in\bN$.

\subsection{Around Hirohisma's theorem}

A straightforward application of Theorem \ref{thm-direct} is the following beautiful result first proved by Hiroshima in 2003  (see Section 4 for the complete form of Hiroshima's theorem and is relevance in quantum physics). 

\vskip 10pt
\begin{cor}\label{Hiroshima} Let $H=[A_{s,t}]\in \bM_{\alpha n}^+$ be written in $\alpha\times \alpha$  Hermitian blocks in $\bM_n$ and let $\Delta=\sum_{s=1}^{\alpha}{A_{s,s}}$ be its partial trace.   Then,  we have
\begin{equation*}
\left\|H \right\|\le \left\| \Delta\right\|
\end{equation*} for all symmetric norms.
\end{cor}

\vskip 10pt
\begin{proof} By completing $H$ with some zero rows and columns, we may assume that $\alpha=\beta$ is dyadic. Theorem \ref{thm-direct} then implies, with $m=2^{\beta}$,
$$
\| \oplus^m H \| \le \| \oplus^m \Delta \| 
$$
for all symmetric norms, which is equivalent to the claim of the corollary.
\end{proof}

\vskip 10pt
By the Ky Fan principle, Corollary \ref{Hiroshima} is equivalent to the majorisation relation
$$
\sum_{i=1}^j\lambda_{i}(H) \le \sum_{i=1}^j\lambda_{i}(\Delta)
$$
for all $j=1,\ldots,\alpha n$. (we set $\lambda_j(A)=0$ when $A\in\bM_d^+$ and $j>d$).

Theorem \ref{thm-direct} says much more than this majorisation. For instance, we may completes these eigenvalue relations with the following ones.

\vskip 10pt
\begin{cor}\label{cor-eigenvalue1} Let $H=[A_{s,t}]\in \bM_{\alpha n}^+$ be written in $\alpha\times \alpha$  Hermitian blocks in $\bM_n$ and let $\Delta=\sum_{s=1}^{\alpha}{A_{s,s}}$ be its partial trace.   Then,  we have
\begin{equation*}
\lambda_{1+\beta k}(H) \le \lambda_{1+k}( \Delta)
\end{equation*} 
 for all $k=0,\ldots,n-1$,    where $\beta$ is the smallest dyadic number such that $\alpha \le\beta$.
\end{cor}

\vskip 10pt
\begin{proof} By completing $H\in\bM_{\alpha n}^+$ with some zero blocks, we may assume that $H\in\bM_{\beta n}^+$. Theorem \ref{thm-direct} then yields the decomposition
\begin{equation*}
\oplus^{m} H =\frac{1}{\beta} \sum_{k=1}^{\beta} V_k\left(\oplus^{m}\Delta\right) V_k^*
\end{equation*} 
where   $\{V_k\}_{k=1}^{\beta}$  is a family of isometries in $\bM_{m\beta n,mn}$ and $m=2^{\beta}$.
We recall a simple fact, Weyl's theorem: if $Y, Z \in\bM_d$ are Hermitian,  then
\[\lambda_{r+s+1}(Y+Z)\le \lambda_{r+1}(Y) + \lambda_{s+1}(Z)\]
for all nonnegative integers $r,s$ such that $r + s\le d-1$. When $Y,Z$ are positive, this still holds for all nonnegative integers $r,s$ with our convention ($\lambda_j(A)=0$ when $A\in\bM_d^+$ and $j>d$). From the previous decomposition we thus infer
\begin{equation}\label{F-1}
\lambda_{1+\beta k}\left(\oplus^m H\right) \le \lambda_{1+ k}\left(\oplus^m \Delta\right)
\end{equation}
for all $k=0,1,\ldots$. Then,  observe that for all $A\in\bM_d^+$ and all $j=0,1,\ldots$,
\begin{equation}\label{F-2}
\lambda_{1+j}\left(\oplus^m A\right) = \lambda_{\langle(1+j)/m\rangle}(A)
\end{equation}
where $\langle u\rangle$ stands for the smallest integer greater than or equal to $u$. Combining \eqref{F-1} and \eqref{F-2} 
 we get
\begin{equation*}
\lambda_{\langle (1+\beta j)/m\rangle}\left(H\right) \le \lambda_{\langle (1+ j)/m\rangle}\left( \Delta\right)
\end{equation*}
for all $j=0,1,\ldots$. Taking $j=km$, $k=0,1\ldots$ completes the proof.
\end{proof}

\vskip 10pt
The above proof actually shows more eigenvalue inequalities.

\vskip 10pt
\begin{cor}\label{cor-eigenvalueC} Let $H=[A_{s,t}]\in \bM_{\alpha n}^+$ be written in $\alpha\times \alpha$  Hermitian blocks in $\bM_n$ and let $\Delta=\sum_{s=1}^{\alpha}{A_{s,s}}$ be its partial trace.   Then,  we have
$$
\lambda_{1+\beta k}( S) \le \frac{1}{\beta}\left\{ \lambda_{1+k_1}\left( \Delta \right)
+\cdots+ \lambda_{1+k_{\beta}}\left( \Delta \right)\right\}
$$
where $k_1+\cdots+k_{\beta}=\beta k$ and $\beta$ is the smallest dyadic number such that $\alpha \le\beta$.
\end{cor}

\vskip 10pt
Corollary \ref{Hiroshima} implies the following rearrangement inequality.

\vskip 10pt
\begin{cor}\label{cor-rearr}
 Let $\{S_i\}_{i=1}^{\alpha}$ be a commuting family of Hermitian operators in $\bM_n$ and let $T\in \bM_n^+$. Then,
\begin{equation*}
\left\| \sum_{i=1}^{\alpha} S_iT^2S_i \right\| \le \left\|  \sum_{i=1}^{\alpha} TS^2_iT \right\|
\end{equation*}
for all symmetric norms.
\end{cor}

\vskip 10pt
\begin{proof}
Define a matrix $Z\in\bM_{\alpha n}$ by
$$ Z= XX^*=
\begin{bmatrix} TS_1\\ \vdots \\  TS_{\alpha}
\end{bmatrix}
\begin{bmatrix} S_1T & \cdots  &S_{\alpha}T
\end{bmatrix}.
$$
Hence $Z=[TS_iS_jT]$ is positive and partitioned in Hermitian blocks in $\bM_n$, with diagonal blocks $TS_i^2T$, $1\le i\le \alpha$. Thus, for all symmetric norms,
$$
\| Z\| \le \left\|\sum_{i=1}^{\alpha} TS_i^2T \right\|
$$
Since $XX^*$ and $X^*X$ have  same symmetric norms for any rectangular matrix $X$, we infer
\begin{equation*}
\left\|\sum_{i=1}^{\alpha}  S_iT^2S_i \right\| \le \left\|\sum_{i=1}^{\alpha} TS_i^2T \right\|
\end{equation*}
as claimed.
\end{proof}

\vskip 10pt
From Corollary \ref{cor-eigenvalue1} we similarly get the next one.

\vskip 10pt
\begin{cor}\label{cor-rearr2}
 Let $\{S_i\}_{i=1}^{\alpha}$ be a commuting family of Hermitian operators in $\bM_n$ and let $T\in \bM_n^+$. Then,
\begin{equation*}
\lambda_{1+\beta k}\left( \sum_{i=1}^{\alpha} S_iT^2S_i \right) \le \lambda_{1+k} \left(  \sum_{i=1}^{\alpha} TS^2_iT \right)
\end{equation*}
 for all $k=0,\ldots,n-1$,    where $\beta$ is the smallest dyadic number such that $\alpha \le\beta$.
\end{cor}

\subsection{Around Rotfel'd inequality}

\noindent
 Given two Hermitian matrices $A$, $B$ in $\bM_n$ and a concave function $f(t)$ defined on the real line,
\begin{equation}\label{VN}
{\mathrm{Tr\,}} f\left(\frac{ A+B}{2}\right) \ge {\mathrm{Tr\,}}\frac{f(A) +f(B)}{2}
\end{equation}
and, if further $f(0)\ge 0$ and both $A$ and $B$ are positive semi-definite,
\begin{equation}\label{Rot}
{\mathrm{Tr\,}} f(A+B) \le {\mathrm{Tr\,}} f(A)+ {\mathrm{Tr\,}}f(B).
\end{equation}
The first inequality goes back to von-Neumann in the 1920's, the second is more subtle and has been proved only in 1969 by Rotfel'd \cite{Rot}. These trace inequalities are matrix versions of  obvious scalar inequalities. Theorem \ref{thm-direct} yields a refinement of the Rotfel'd inequality for  families of positive operators $\{A_i\}_{i=1}^{\alpha}$ by considering these operators as the diagonal blocks of a partitioned matrix $H$ as follows.

\vskip 10pt
\begin{cor}\label{cor-rot} Let $H=[A_{s,t}]\in \bM_{\alpha n}^+$ be written in $\alpha\times \alpha$  Hermitian blocks in $\bM_n$. Then, we have 
$$
{\mathrm{Tr}}\, f\left(\sum_{s=1}^{\alpha} A_{s,s}\right)
  \le  {\mathrm{Tr}}\,f(H)
 \le  \sum_{s=1}^{\alpha}{\mathrm{Tr}}\,f(A_{s,s})
$$
for all concave functions $f(t)$ on $\bR^+$ such that $f(0)\ge 0$.
\end{cor}

\vskip 10pt
\begin{proof}
Note that if these inequalities hold for a non-negative concave function $f(t)$ with $f(0)=0$, then they also hold for the function $f(t)+c$ for any constant $c>0$. Therefore it suffice to consider  concave functions vanishing at the origine. This assumption entails that
\begin{equation}\label{FF-1}
f(VAV^*)=Vf(A)V^*
\end{equation}
for all $A\in\bM_n^+$ and all isometries $V\in\bM_{m,n}$. By Lemma \ref{BL-lemma} we have a decomposition
$$
H=\sum_{s=1}^{\alpha} V_sA_{s,s}V_s^*
$$
for some isometries $V_s\in\bM_{\alpha n,n}$. Inequality \eqref{Rot} then yields
$$
 {\mathrm{Tr}}\,f(H) \le \sum_{s=1}^{\alpha} {\mathrm{Tr}}\,f(V_sA_{s,s}V_s^*)
$$
and using \eqref{FF-1} establishes the second inequality. To prove the first inequality, we use Theorem \ref{thm-direct}. By completing $H$ with some zero blocks (we still suppose $f(0)=0$) we may assume that $\alpha=\beta$ is dyadic. Thus we have a decomposition
\begin{equation*}
\oplus^{m} H =\frac{1}{\beta} \sum_{k=1}^{\beta} V_k\left(\oplus^{m}\Delta\right) V_k^*
\end{equation*} 
where   $\{V_k\}_{k=1}^{\beta}$  is a family of isometries in $\bM_{m\beta n,mn}$ and $m=2^{\beta}$. Inequality \eqref{VN} then gives
$$ {\mathrm{Tr}}\, f\left(\oplus^{m} H \right) \ge 
\frac{1}{\beta}\sum_{k=1}^{\beta}{\mathrm{Tr}}\,f\left( V_k\left(\oplus^{m}\Delta\right) V_k^*\right)
$$
and using \eqref{FF-1} we obtain
$$ {\mathrm{Tr}}\ f\left(\oplus^{m} H \right) \ge {\mathrm{Tr}}\ f\left(\oplus^{m} \Delta \right).
$$
The proof is completed by dividing both sides by $m$.
\end{proof}

\vskip 10pt
\begin{remark} The first inequality of Corollary \ref{cor-rot} is actually equivalent to Corollary \ref{Hiroshima}
by a well-known majorisation  principle  for convex/concave functions. The above proof does not require this principle. The simplest case of Corollary \ref{cor-rot} is the double inequality
$$
 f(a_1+\cdots+a_n) \le {\mathrm{Tr}}\, f(A)
\le f(a_1)+\cdots +f(a_n)
$$
for all $A\in\bM_n^+$ with diagonal entries $a_1,\ldots,a_n$.

\end{remark}

\vskip 10pt
\begin{remark} A special case of Corollary \ref{cor-rot} refines a well-known determinantal inequality. Taking as a concave function on $\bR^+$, $f(t)=\log(1+t)$, we obtain:
Let $A,B\in \bM_n^+$. Then, for any Hermitian $X\in\bM_n$  such that 
$$H=\begin{bmatrix} A&X \\X&B\end{bmatrix}$$
is positive semi-definite, we have 
\begin{equation*}
\det(I+A+B)
\le \det(I+H) \le \det(I+A)\det(I+B).
\end{equation*}
This was noted in \cite{BLL2}.
\end{remark}

\section{Proof of Theorem 2.1}

\noindent
A Clifford algebra ${\mathcal{C}}_{\beta}$ is the associative real algebra generated by $\beta$ elements $q_1,\ldots,q_{\beta}$ satisfying   the canonical anticommutation relations $q_i^2=1$ and
$$
q_iq_j +q_jq_i =0 
$$
for $i\neq j$. This structure was introduced by Clifford in \cite{Cli}. It turned out to be of great importance in quantum theory and operator algebras, for instance see  the survey \cite{Der}. From the relation
$$
\begin{pmatrix}
0&1 \\
1&0
\end{pmatrix}
\begin{pmatrix}
1&0 \\
0&-1
\end{pmatrix}
+
\begin{pmatrix}
1&0 \\
0&-1
\end{pmatrix}
\begin{pmatrix}
0&1 \\
1&0
\end{pmatrix}=0
$$
we infer a representation of  ${\mathcal{C}}_{\beta}$ as a a real subalgebra of $M_{2^{\beta}}=\otimes^{\beta} \bM_2$ by mapping the generators $q_j\mapsto Q_j$, $1\le j\le \beta$, where
\begin{equation}\label{cli-gen}
Q_j=\left\{\otimes^{j-1}\begin{pmatrix}
1&0 \\
0&-1
\end{pmatrix}\right\}\otimes
\begin{pmatrix}
0&1 \\
1&0
\end{pmatrix}
\otimes
\left\{\otimes^{\beta-j}\begin{pmatrix}
1&0 \\
0&1
\end{pmatrix}
\right\}.
\end{equation}
We use these matrices in the following proof of Theorem \ref{thm-direct}.

\vskip 10pt
\begin{proof} 
First, replace the positive block matrix $H=[A_{s,t}]$ where $1\le s,t,\le \beta$
and all blocks are Hermitian by a bigger one in which each block in counted $2^{\beta}$ times :
$$G= [G_{s,t}]:= \left[I_{2^{\beta}}\otimes A_{s,t}\right]= \left[\oplus^{2^{\beta}} A_{s,t}\right]$$
where $I_r$ stands for the identity of $\bM_r$.
Thus $G\in\bM_{\beta 2^{\beta}n}$ is written in $\beta$-by-$\beta$ blocks in $\bM_{2^{\beta}n}$. Then perform a unitary congruence with the unitary $W\in\bM_{\beta2^{\beta}n}$ defined as
\begin{equation}\label{unitaryW}
W=\bigoplus_{j=1}^{\beta} \left\{Q_j \otimes I_n\right\}
\end{equation}
where $Q_j$ is given by \eqref{cli-gen}, $1\le j\le \beta$.
Thanks to the anticommutation relation for each pair of summands in
 \eqref{unitaryW},
$$
\left\{Q_j \otimes I_n\right\}\left\{Q_l \otimes I_n\right\}+\left\{Q_l \otimes I_n\right\}\left\{Q_j \otimes I_n\right\}=0, \quad j\neq l,
$$
 the block matrix (with $W=W^*$)
\begin{equation}\label{omega1}
\Omega:=WGW^*=[ \Omega_{s,t}]
\end{equation}
satisfies the following  :
For $1\le s<t\le \beta$,
\begin{equation}\label{omega2}
\Omega_{s,t}=-\Omega_{t,s}.
\end{equation}

Next, consider the reflexion matrix
$$
J_1=\frac{1}{\sqrt{2}} \begin{pmatrix} 1&1 \\1&-1 \end{pmatrix}
$$
and define inductively for all integers $p>1$, a reflexion
$$
J_p=\frac{1}{\sqrt{2}} \begin{pmatrix} J_{p-1} & J_{p-1}  \\ J_{p-1} &- J_{p-1}  \end{pmatrix},
$$
that is $J_p=\otimes^p J_1$. Observe, that given any  matrix $S\in\bM_{2^p}$, $S=[s_{i,j}]$, such that $s_{i,j}=-s_{i,j}$ for all $i\neq j$, the matrix
$$
T=J_pSJ_p^*
$$
has its diagonal entries $t_{j,j}$ all equal to the normalized trace
$2^{-p}{\mathrm{Tr\,}} S$. Indeed, letting $J_p=[z_{i,j}]$, 
\begin{align*}
t_{j,j} &=\sum_{k} z_{j,k} \left(\sum_{l}s_{k,l}z_{l,j}\right) \\
&=\sum_{k,l} z_{j,k} s_{k,l} z_{l,j} \\
&=\sum_{k} z_{j,k} s_{k,k} z_{k,j} +\sum_{k\neq l}  z_{j,k} s_{k,l} z_{l,j}\\
&=2^{-p}{\mathrm{Tr\,}} S +\sum_{k<l}\left( s_{k,l} z_{j,k}z_{l,j}+s_{l,k} z_{l,j}z_{k,j}\right)\\
&=2^{-p}{\mathrm{Tr\,}} S.
\end{align*}

Now, since we assume that $\beta=2^p$ for some integer $p$, we may  perform a unitary congruence to the matrix $\Omega$ in \eqref{omega1} with the unitary matrix
$$
R_p=J_p\otimes I_{2^{\beta}}\otimes I_n
$$
and, making use of \eqref{omega2} and the above property of $J_p$, we note that 
$
R_p\Omega R_p^*
$
has  its $\beta$ diagonal blocks $(R_p\Omega R_p^*)_{j,j}$, $1\le j\le \beta$, all equal to the matrix $D\in\bM_{2^{\beta}n}$,
$$
D =\frac{1}{\beta}\sum_{s=1}^{\beta} \left\{ \oplus^{2^{\beta}}A_{s,s}\right\}.
$$
Thanks to the decomposition  of Lemma \ref{BL-lemma},  there exist some isometries $U_k\in\bM_{\beta2^{\beta}n, 2^{\beta}n}$, $1\le k\le \beta$, such that
$$
\Omega=\sum_ {k=1}^{\beta} U_k D U_k^*.
$$
 Since $\Omega$ is unitarily equivalent to $\oplus^{2^{\beta}} H$, that is $\Omega=V^*(\oplus^{2^{\beta}} H)V$ for some unitary $V\in \bM_{\beta 2^{\beta}n, 2^{\beta}n}$, we get
$$
\oplus^{2^{\beta}} H = \sum_ {k=1}^{\beta} VU_k D U_k^*V^*
$$
wich is the claim of Theorem \ref{thm-direct} by setting $VU_k=:V_k$, $1\le k\le \beta$, as $2^{\beta}=m$, and $D=\frac{1}{\beta}\oplus^m\Delta$.
\end{proof}

\section{Comments}

\subsection{Complex matrices and small partitions}

If one uses isometries with complex entries, then, in case of partitions into a small number of $\beta\times\beta$ blocks, the number $m$ of copies in the direct sum $\oplus^m H$ and $\oplus^m \Delta$ can be reduced. For $\beta=2$, Theorem \ref{thm-four} shows that it suffices to take $m=1$. For $\beta=3$ or $\beta=4$ the following result holds \cite{BLL2}.

\vskip 10pt
\begin{theorem}\label{thm-quaternion} Let $H=[A_{s,t}]\in \bM_{\beta n}^+(\bC)$ be written in  Hermitian blocks in $\bM_n(\bC)$  with $\beta\in\{3,4\}$ and let $\Delta=\sum_{s=1}^{\beta}A_{s,s}$ be its partial trace. Then, 
\begin{equation*}
H\oplus H =\frac{1}{4}\sum_{k=1}^4 V_k\left(\Delta\oplus\Delta\right)V_k^*
\end{equation*} for some isometries $V_k\in\bM_{2\beta n,2n}(\bC)$, $k=1,2,3,4$.
\end{theorem}

\vskip 10pt
Likewise for Theorem \ref{thm-four}, we must consider isometries with complex entries, even for a full matrix $H$ with real entries. The proof makes use of quaternions and thus  confines to $\beta\le 4$.

\subsection{Separability criterion}

\vskip10pt
Let ${\mathcal{H}}$ and ${\mathcal{F}}$ be two  finite dimensional Hilbert spaces that may be either  real spaces,   identified to $\bR^n$ and $\bR^m$, or complex spaces,  identified to $\bC^n$ and $\bC^m$. The space of operators on ${\mathcal{H}}$, denoted by ${\mathrm{B}}({\mathcal{H}})$, is  identified with the matrix algebra $\bM_n$ (with real or complex entries according the nature of ${\mathcal{H}}$). A positive (semi-definite) operator $Z$ on the tensor product space ${\mathcal{H}}\otimes{\mathcal{F}}$ is said to be separable if it can be decomposed as a sum of tensor products of positive operators,
\begin{equation}\label{sum-tens2}
Z=\sum_{j=1}^k A_j\otimes B_j
\end{equation}
where $A_j$'s are positive operators on ${\mathcal{H}}$ and so  $B_j$'s are on ${\mathcal{F}}$. 
It is  difficult in general to determine if a given
positive operator in the  matrix algebra $\bM_n\otimes\bM_m$
  is separable or not, though some theoretical criteria do exist \cite{Hor}, \cite{Chen-Wu}. The  partial trace  of  $Z$ with respect to ${\mathcal{H}}$ is  the operator acting on ${\mathcal{F}}$,
$$
{\mathrm{Tr}}_{\mathcal{H}} Z= \sum_{j=1}^k ({\mathrm{Tr}}A_j)B_j.
$$
These notions have their own mathematical interest and moreover play a fundamental role in the description of bipartite systems in quantum theory, see \cite[Chapter 10]{Petz}, where the positive operators act on complex spaces and are usually normalized with trace one and  called states. Thus a separable state is an operator of the type \eqref{sum-tens2} with ${\mathrm{Tr\,}} Z=1$. The richness of the mathematical theory of separable operators/states and their application in quantum physics is apparent in many places in the literature, for instance in \cite{Hor} and \cite{AlSh}. Nielsen and Kempe in 2001 proved a majorisation separability criterion \cite{NK}. It can be stated as the following norm comparison.

\vskip 10pt
\begin{theorem}\label{part-trace} Let  $Z$ be a separable state on the tensor product   of two  finite dimensional  Hilbert spaces ${\mathcal{H}}$ and ${\mathcal{F}}$. Then, for all symmetric norms,
$$
\| Z\| \le \left\| {\mathrm{Tr}}_{\mathcal{H}} Z \right\|.
$$
\end{theorem}

\vskip 10pt
Regarding ${\mathrm{B}}({\mathcal{H}}\otimes{\mathcal{F}})$ as $\bM_n(\bM_m)$, an operator $Z\in{\mathrm{B}}({\mathcal{H}}\otimes{\mathcal{F}})$ is written as a block-matrix $Z=[Z_{i,j}]$ with $Z_{i,j}\in\bM_m$, $1\le i,j\le n$. The partial trace of $Z$  with respect to ${\mathcal{H}}$ is then the sum of the diagonal blocks,
$$
{\mathrm{Tr}}_{\mathcal{H}} Z=\sum_{j=1}^n Z_{j,j}.
$$
This observation makes obvious that  Theorem \ref{part-trace} is a straightforward consequence of  Corollary \ref{Hiroshima} whenever the factor ${\mathcal H}$ is a {\bf real} Hilbert space. Indeed, we then have
$$
Z=\sum_{j=1}^k A_j\otimes B_j
$$
where, for each index $j$,  $A_j\in\bM_n(\bR)$  and  $B_j$ is Hermitian in $\bM_m$,  so that $A_j\otimes B_j$ can be regarded as an element of $\bM_n( \bM_m)$ formed of Hermitian blocks.

From Corollary \ref{cor-eigenvalue1}, we may complete the majorisation of Theorem \ref{part-trace}, when a factor is a real space with a few more eigenvalue estimate as stated in the next corollary.

\vskip 10pt
\begin{cor} Let  $Z$ be a separable positive operator on the tensor product of two finite dimensional  Hilbert space ${\mathcal{H}}\otimes{\mathcal{F}}$ with a {\bf real} factor ${\mathcal{H}}$. Then,
$$
\lambda_{1+\beta k}( Z) \le \lambda_{1+k}\left( {\mathrm{Tr}}_{\mathcal{H}} Z \right)
$$
for all $k=0,\ldots,\dim{\mathcal{F}}-1$,    where $\beta$ is the smallest dyadic number such that $\dim{\mathcal{H}} \le\beta$.
\end{cor}

\vskip 10pt
Similarly, from Corollary \ref{cor-eigenvalueC}, we actually have a larger set of eigenvalue inequalities.

\vskip 10pt
\begin{cor} Let  $Z$ be a separable positive operator on the tensor product of two finite dimensional  Hilbert space ${\mathcal{H}}\otimes{\mathcal{F}}$ with a {\bf real} factor ${\mathcal{H}}$. Then,
$$
\lambda_{1+\beta k}( Z) \le \frac{1}{\beta}\left\{ \lambda_{1+k_1}\left( {\mathrm{Tr}}_{\mathcal{H}} Z  \right)
+\cdots+ \lambda_{1+k_{\beta}}\left( {\mathrm{Tr}}_{\mathcal{H}} Z  \right)\right\}
$$
for all $k=0,\ldots,\dim{\mathcal{F}}-1$,
where $k_1+\cdots+k_{\beta}=\beta k$ and $\beta$ is the smallest dyadic number such that $\dim{\mathcal{H}} \le\beta$.
\end{cor}

\vskip 10pt
Of course, confining to real spaces is a severe restriction. It would be desirable to obtain similar estimates for usual complex spaces. A related problem would be to obtain a decomposition like in Theorem \ref{thm-direct} for the class of partitioned matrices considered in the full form of Hiroshima's theorem: 
\vskip 10pt\noindent
{\it If $A=[A_{s,t}]$ is a positive matrix partitioned in $\alpha\times \alpha$ blocks such that $B=[B_{s,t}]:=[A_{s,t}^T]$ is positive too, then the majorisation of Corollary \ref{Hiroshima} holds.}
\vskip 10pt\noindent
 Here $X^T$ means the transposed matrix. This statement extends  Corollary \ref{Hiroshima} and implies Theorem \ref{part-trace}.

\vskip 10pt

J.-C. Bourin,

Laboratoire de math\'ematiques,

Universit\'e de Franche-Comt\'e,

25 000 Besan\c{c}on, France.

jcbourin@univ-fcomte.fr

\vskip 10pt
Eun-Young Lee

 Department of mathematics,

Kyungpook National University,

 Daegu 702-701, Korea.

eylee89@ knu.ac.kr


\begin{thebibliography}{99}
{\small

\bibitem{AlSh} E.\ Alfsen and F.\ Shultz,
Unique decompositions, faces, and automorphisms
of separable states, {\it J.\  Math.\ Phys.\  }\textbf{51}, 052201 (2010)




\bibitem{Bhatia} R.\ Bhatia, Positive Definite Matrices, 




\bibitem{BH1} J.-C.\ Bourin and F.\ Hiai,
Norm and anti-norm inequalities for positive semi-definite matrices,
{\it Internat. J. Math.} \textbf{22} (2011), 1121-1138.



 \bibitem{BL1}
J.-C. Bourin and E.-Y. Lee, Unitary orbits of Hermitian operators with convex and
concave functions, {\it Bull.\ London Math.\ Soc.}, in press.

\bibitem{BLL1}
J.-C. Bourin, E.-Y. Lee and  M.\ Lin, On a decomposition lemma for positive semi-definite block-matrices, {\it Linear Algebra Appl.}  \textbf{437} (2012), 1906-1912.

\bibitem{BLL2} J.-C. Bourin, E.-Y. Lee and  M.\ Lin, Positive matrices partitioned into a small number of
Hermitian blocks, preprint.

\bibitem{Chen-Wu}
K.\ Chen and L.-A.\ Wu, A matrix realignment method for recognizing entanglement, {\it Quantum Inf.\ Comput.}\ \textbf{3}  (2003), 193-202.

\bibitem{Cli} Clifford, Applications of Grassmann' extensive algebra, {\it Amer.\ Journ.\ Math.}\ \textbf{1} (1878), 350-358.

\bibitem{Der} J.\ Derezi\'nski, {\it Introduction to representations of the canonical commutation and anticommutation relations}, Lect.\ Note Phys.\ {\bf 695}, 65-145 (2006), Springer.


\bibitem{H} T.\ Hiroshima, Majorization Criterion for Distillability of a Bipartite Quantum State, (ArXiv: quantum-ph, 2003).


\bibitem{Hor} M.\ Horodecki, P.\ Horodecki, R.\ Horodecki, Separability of mixed states: necessary and sufficient conditions,
{\it Phys.\ Lett.\ A} \textbf{223} (1996) l-8.

\bibitem{NK} M.\ A.\ Nielsen and J.\ Kempe,  Separable states are more disordered globally than locally, {\it Phys.\ Rev.\ Lett.}\ \textbf{86} (2001) 5184-5187.

\bibitem{Petz} D.\ Petz, Matrix Analysis with some applications, $<$http://www.math.hu/petz$>$.


\bibitem{Rot}  S. Ju. Rotfel'd, The singular values of a
sum of completely continuous operators,  \textit{ Topics in
Mathematical Physics, Consultants Bureau}, Vol. \textbf{3} (1969)
73-78.


}

\end{thebibliography}
\end{document}